\documentclass[psamsfonts,reqno]{amsart}
\usepackage{amssymb,eucal,graphics,latexsym}

\newcommand{\SUB}{\mathbf{SUB}}

\newcommand{\upw}{\mathbin{\uparrow}}

\newcommand{\bP}{\mathbb{P}}
\newcommand{\bD}{\mathbb{D}}
\newcommand{\PP}{\mathcal{P}}
\newcommand{\jz}{$\langle\vee,0\rangle$}
\newcommand{\jzs}{\jz-sem\-i\-lat\-tice}

\newcommand{\jh}{join-ho\-mo\-mor\-phism}
\newcommand{\mh}{meet-ho\-mo\-mor\-phism}

\newcommand{\jirr}{join-ir\-re\-duc\-i\-ble}
\newcommand{\mirr}{meet-ir\-re\-duc\-i\-ble}
\newcommand{\jirry}{join-ir\-re\-duc\-i\-bil\-i\-ty}

\newcommand{\CD}{\mathcal{D}}
\newcommand{\DD}{\mathbin{D}}
\newcommand{\rd}[1]{[{#1}]^{\DD}}
\newcommand{\pup}[1]{\textup{(}{#1}\textup{)}}

\newcommand{\ol}[1]{\overline{#1}}

\newcommand{\fsi}[1]{\{1,\dots,#1\}}
\newcommand{\fso}[1]{\{0,\dots,#1\}}

\newcommand{\es}{\varnothing}

\newcommand{\tr}{\vartriangleleft}

\newcommand{\utr}{\trianglelefteq}
\newcommand{\nutr}{\ntrianglelefteq}
\newcommand{\dtr}{\mathbin{\vartriangleleft\kern-10pt {\lower
3pt\hbox{$\scriptscriptstyle\neq$}}\kern3pt}}

\newcommand{\set}[1]{\{{#1}\}}
\newcommand{\setm}[2]{\set{{#1}\mid{#2}}}
\newcommand{\seq}[1]{\langle#1\rangle}
\newcommand{\famm}[2]{\seq{{#1}\mid{#2}}}

\newcommand{\St}{\textup{(S)}}
\newcommand{\Ud}{\textup{(U)}}
\newcommand{\Bo}{\textup{(B)}}

\newcommand{\Ht}[1]{\textup{(H$_{#1}$)}}
\newcommand{\Lt}{\textup{(L$_{2}$)}}

\DeclareMathOperator{\lh}{length}

\newcommand{\Stj}{\textup{(S$_{\mathrm{j}}$)}}
\newcommand{\Udj}{\textup{(U$_{\mathrm{j}}$)}}
\newcommand{\Boj}{\textup{(B$_{\mathrm{j}}$)}}

\newcommand{\fr}{\mathbf{r}}

\newcommand{\hL}{{\widehat{L}}}

\newcommand{\into}{\hookrightarrow}
\newcommand{\onto}{\twoheadrightarrow}

\newcommand{\eps}{\varepsilon}

\DeclareMathOperator{\J}{J}
\DeclareMathOperator{\Fil}{Fil}

\newcommand{\Co}{\mathbf{Co}}

\numberwithin{equation}{section}

\theoremstyle{plain}
\newtheorem{lemma}{Lemma}[section]
\newtheorem{theorem}[lemma]{Theorem}
\newtheorem{proposition}[lemma]{Proposition}
\newtheorem{corollary}[lemma]{Corollary}
\newtheorem{example}[lemma]{Example}
\newtheorem{claim}{Claim}

\newtheorem*{stat}{\name}
\newcommand{\name}{testing}

\theoremstyle{definition}
\newtheorem{definition}[lemma]{Definition}
\newtheorem{notation}[lemma]{Notation}
\newtheorem{problem}{Problem}

\theoremstyle{remark}
\newtheorem{remark}[lemma]{Remark}
\newtheorem*{note}{Note}

\newenvironment{all}[1]{\renewcommand{\name}{#1}\begin{stat}}
                        {\end{stat}}

\newcommand{\qedc}{{\qed}~{\rm Claim~{\theclaim}.}}
\newcommand{\qedsc}{{\qed}~{\rm Claim.}}

\newenvironment{cproof}
{\begin{proof}[Proof of Claim.]}
{\qedc\renewcommand{\qed}{}\end{proof}}

\begin{document}

\title[Lattices of order-convex sets, II]%
{Sublattices of lattices of order-convex sets, II.\\
Posets of finite length}

\author[M.~Semenova]{Marina Semenova}
\address[M.~Semenova]{Institute of Mathematics of
the Siberian Branch of RAS\\
Acad. Koptyug prosp. 4\\
630090 Novosibirsk\\
Russia}
\email{semenova@math.nsc.ru}

\author[F.~Wehrung]{Friedrich Wehrung}
\address[F.~Wehrung]{CNRS, UMR 6139\\
D\'epartement de Math\'ematiques\\
Universit\'e de Caen\\
14032 Caen Cedex\\
France}
\email{wehrung@math.unicaen.fr}
\urladdr{http://www.math.unicaen.fr/\~{}wehrung}
\date{\today}


\thanks{The first author was supported by INTAS grant no.
YSF: 2001/1-65, by RFBR grants no. 99-01-00485 and 01-01-06178,
by GA UK grant no.~162/1999, and by GA CR grant no. 201/99.
The second author was partially supported by the
Fund of Mobility of the Charles University (Prague), by FRVS grant
no. 2125, by institutional grant CEZ:J13/98:113200007, and by the
Barrande program.}

\begin{abstract}
For a positive integer $n$, we denote by $\SUB$ (resp., $\SUB_n$)
the class of all lattices that can be embedded into the lattice
$\Co(P)$ of all order-convex subsets of a partially ordered set $P$
(resp., $P$ of length at most $n$). We prove the following results:
\begin{itemize}
\item[(1)] $\SUB_n$ is a finitely based variety, for any
$n\geq1$.

\item[(2)] $\SUB_2$ is locally finite.

\item[(3)] A finite atomistic lattice $L$ without $\DD$-cycles
belongs to $\SUB$ if{f} it belongs to $\SUB_2$; this result does not
extend to the nonatomistic case.

\item[(4)] $\SUB_n$ is not locally finite for $n\geq3$.

\end{itemize}

\end{abstract}

\maketitle

\section{Introduction}\label{S:Intro}

For a partially ordered set (from now on \emph{poset}) $(P,\utr)$, a
subset $X$ of $P$ is \emph{order-convex}, if $x\utr z\utr y$ and
$\set{x,y}\subseteq X$ implies that $z\in X$, for all $x$, $y$,
$z\in P$. The set $\Co(P)$ of all order-convex subsets of $P$ forms a
lattice under inclusion. It gives an important example of
\emph{convex geometry}, see K.\,V. Adaricheva, V.\,A. Gorbunov, and
V.\,I. Tumanov \cite{AGT}. In M. Semenova and F. Wehrung
\cite{SeWe1}, the following result is proved:

\begin{all}{Theorem}
The class $\SUB$ of all lattices that can be embedded into
some $\Co(P)$ is a variety.
\end{all}

This implies the nontrivial result that \emph{every homomorphic image
of a member of $\SUB$ belongs to $\SUB$}. It is in fact proved in
\cite{SeWe1} that the variety $\SUB$ is \emph{finitely based}, it is
defined by three identities that are denoted by \St, \Ud, and \Bo.

In the present paper, we extend this result to the class $\SUB_n$ of
all lattices that can be embedded into $\Co(P)$ for some poset $P$
of length $n$, for a given positive integer $n$:

\begin{all}{Theorem~\ref{T:SUBnVar}}
The class $\SUB_n$ is a finitely based variety, for every positive
integer~$n$.
\end{all}

It is well-known that for $n=1$, the class $\SUB_n$ is the variety of
all \emph{distributive} lattices. This fact is contained in
G. Birkhoff and M.\,K. Bennett \cite{BB}.

For $n=2$, $\SUB_n=\SUB_2$ is much more interesting, it is the
variety of all lattices that can be embedded into some $\Co(P)$
\emph{without $\DD$-cycle on its atoms}. We find a simple finite
set of identities characterizing $\SUB_2$, see Theorem~\ref{T:SUB2}.
In addition, we prove the following results:

\begin{itemize}\em
\item[---] The variety $\SUB_2$ is locally finite \pup{see
Theorem~\textup{\ref{T:SUB2lf}}}, and we provide an explicit upper
bound for the cardinality of the free lattice on $m$ generators in
$\SUB_2$.

\item[---] A finite atomistic lattice without $\DD$-cycle belongs to
$\SUB$ if{f} it belongs to $\SUB_2$ \pup{see
Proposition~\textup{\ref{P:AtSUB2}}}.
\end{itemize}

We also prove that $\SUB_n$ is not locally finite for $n\geq3$ (see
Theorem~\ref{T:NonLF3}), and that $\SUB_n$ is a proper subvariety of
$\SUB_{n+1}$ for every $n$ (see Corollary~\ref{C:chainsHn}).

\section{Basic concepts}\label{S:Basic}

We recall some of the definitions and concepts used in \cite{SeWe1}.
For elements $a$, $b$, $c$ of a lattice $L$ such that
$a\leq b\vee c$, we say that the (formal) inequality $a\leq b\vee c$
is a \emph{nontrivial join-cover}, if $a\nleq b,c$. We say that it is
\emph{minimal in $b$}, if $a\nleq x\vee c$ holds, for all $x<b$, and
we say that it is a \emph{minimal nontrivial join-cover}, if it is a
nontrivial join-cover and it is minimal in both $b$ and $c$.

The \emph{join-dependency} relation $\DD=\DD_L$ (see R. Freese, J.
Je\v{z}ek, and J.\,B. Nation~\cite{FJN}) is defined on the
set $\J(L)$ of all \jirr\ elements of $L$ by putting
 \begin{equation}\label{Eq:DefD}
 p\DD q,\text{ if }p\neq q\text{ and }\exists x\text{ such that }
 p\leq q\vee x\text{ holds and is minimal in }q.
 \end{equation}
It is important to observe that $p\DD q$ implies that $p\nleq q$, for
all $p$, $q\in\J(L)$. Furthermore, $p\nleq x$ in \eqref{Eq:DefD}.

We say that $L$ is \emph{finitely spatial} (resp., \emph{spatial})
if every element of $L$ is a join of \jirr\ (resp., completely \jirr)
elements of $L$. It is well known that every dually algebraic
lattice is lower continuous---see Lemma~2.3 in P. Crawley and R.\,P.
Dilworth~\cite{CrDi}, and spatial (thus finitely spatial)---see Theorem
I.4.22 in G.~Gierz \emph{et al.}~\cite{Comp} or Lemma~1.3.2 in V.\,A.
Gorbunov~\cite{Gorb}.

A lattice $L$ is \emph{dually $2$-distributive}, if it satisfies the
identity
 \begin{equation*}
 a\wedge(x\vee y\vee z)=
 (a\wedge(x\vee y))\vee(a\wedge(x\vee z))\vee(a\wedge(y\vee z)).
 \end{equation*}
A stronger identity is the \emph{Stirlitz identity} \St\ introduced
in \cite{SeWe1}:
 \[
 a\wedge(b'\vee c)=
 (a\wedge b')\vee\bigvee_{i<2}
 \Bigl(a\wedge(b_i\vee c)\wedge
 \bigl((b'\wedge(a\vee b_i))\vee c\bigr)\Bigr),
 \]
where we put $b'=b\wedge(b_0\vee b_1)$. Two other important
identities are the \emph{Udav identity}~\Ud,
 \begin{multline*}
 x\wedge(x_0\vee x_1)\wedge(x_1\vee x_2)\wedge(x_0\vee x_2)\\
 =(x\wedge x_0\wedge(x_1\vee x_2))\vee
 (x\wedge x_1\wedge(x_0\vee x_2))\vee
 (x\wedge x_2\wedge(x_0\vee x_1)),
 \end{multline*}
and the \emph{Bond identity} \Bo,
 \begin{align*}
 x\wedge(a_0\vee a_1)\wedge(b_0\vee b_1)=&\bigvee_{i<2}
 \Bigl(\bigl(x\wedge a_i\wedge(b_0\vee b_1)\bigr)\vee
 \bigl(x\wedge b_i\wedge(a_0\vee a_1)\bigr)\Bigr)\\
 &\vee\bigvee_{i<2}\bigl(x\wedge(a_0\vee a_1)\wedge(b_0\vee b_1)
 \wedge(a_0\vee b_i)\wedge(a_1\vee b_{1-i})\bigr).
 \end{align*}
It is proved in \cite{SeWe1} that \emph{a lattice $L$
belongs to $\SUB$ if{f} it satisfies \St, \Ud, and \Bo}.
Although these identities are quite complicated, they have the
following respective consequences, their so-called
\emph{\jirr\ interpretations}, that can be easily visualized
on the poset $P$ in case $L=\Co(P)$ for a poset $P$:
\begin{itemize}
\item[\Stj:] For all $a$, $b$, $b_0$, $b_1$, $c\in\J(L)$, the
inequalities $a\leq b\vee c$, $b\leq b_0\vee b_1$, and $a\neq b$
imply that either $a\leq\ol{b}\vee c$ for some $\ol{b}<b$ or
$b\leq a\vee b_i$ and $a\leq b_i\vee c$ for some $i<2$.

\item[\Udj:] For all $x$, $x_0$, $x_1$, $x_2\in\J(L)$, the
inequalities
$x\leq x_0\vee x_1,x_0\vee\nobreak x_2,x_1\vee\nobreak x_2$ imply
that either $x\leq x_0$ or $x\leq x_1$ or $x\leq x_2$.

\item[\Boj:] For all $x$, $a_0$, $a_1$, $b_0$, $b_1\in\J(L)$, the
inequalities $x\leq a_0\vee a_1,b_0\vee b_1$ imply that either
$x\leq a_i$ or $x\leq b_i$ for some $i<2$ or $x\leq a_0\vee b_0$,
$a_1\vee b_1$ or $x\leq a_0\vee b_1$, $a_1\vee b_0$.
\end{itemize}

It is proved in \cite{SeWe1} that \St\ implies \Stj, \Ud\ implies
\Udj, and \Bo\ implies \Boj.

A \emph{Stirlitz track} of $L$ is a pair
$(\famm{a_i}{0\leq i\leq n},\famm{a'_i}{1\leq i\leq n})$, where the
$a_i$-s and the $a'_i$-s are \jirr\ elements of $L$ that satisfy the
following relations:
\begin{enumerate}
\item the inequality $a_i\leq a_{i+1}\vee a'_{i+1}$ holds, for all
$i\in\fso{n-1}$, and it is a minimal nontrivial join-cover;

\item the inequality $a_i\leq a'_i\vee a_{i+1}$ holds, for all
$i\in\fsi{n-1}$.
\end{enumerate}

For a poset $P$, the \emph{length of $P$}, denoted by $\lh P$, is
defined as the supremum of the numbers $|C|-1$, where $C$ ranges over
the finite subchains of $P$. We say that $P$ with predecessor
relation $\prec$ is \emph{tree-like}, if it has no infinite bounded
chain and between any points $a$ and $b$ of $P$ there exists at most
one finite sequence $\famm{x_i}{0\leq i\leq n}$ with distinct entries
such that $x_0=a$, $x_n=b$, and either $x_i\prec x_{i+1}$ or
$x_{i+1}\prec x_i$, for all $i\in\fso{n-1}$.

\section{The identity \Lt}\label{S:Letna}

Let \Lt\ be the following lattice-theoretical identity:
 \begin{multline*}
 a\wedge\Bigl(\bigl(b\wedge(c\vee c')\bigr)\vee b'\Bigr)=\\
 \bigl(a\wedge b\wedge(c\vee c')\bigr)\vee
 \Bigl(a\wedge\bigl((b\wedge c)\vee b'\bigr)\Bigr)\vee
 \Bigl(a\wedge\bigl((b\wedge c')\vee b'\bigr)\Bigr).
 \end{multline*}
Taking $b=c\vee c'$ implies immediately the following:

\begin{lemma}\label{L:Lt22DD}
The identity \Lt\ implies dual $2$-distributivity.
\end{lemma}

In order to find an alternative formulation for \Lt\ and many other
identities, it is convenient to introduce the following definition.

\begin{definition}\label{D:GoodGen}
A subset $\Sigma$ of a lattice $L$ is a \emph{join-seed}, if the
following assertions hold:
\begin{enumerate}
\item $\Sigma\subseteq\J(L)$;

\item every element of $L$ is a join of elements of $\Sigma$;

\item for all $p\in\Sigma$ and all $a$, $b\in L$ such that
$p\leq a\vee b$ and $p\nleq a,b$, there are $x\leq a$ and $y\leq b$
both in $\Sigma$ such that $p\leq x\vee y$ is minimal in $x$ and $y$.
\end{enumerate}
\end{definition}

Two important examples of join-seeds are provided by the following:

\begin{lemma}\label{L:JoinSeeds}
Any of the following assumptions implies that the subset $\Sigma$ is
a join-seed of the lattice $L$:
\begin{enumerate}
\item $L=\Co(P)$ and $\Sigma=\setm{\set{p}}{p\in P}$, for some
poset $P$.

\item $L$ is a dually $2$-distributive, complete, lower continuous,
finitely spatial lattice, and $\Sigma=\J(L)$.
\end{enumerate}
\end{lemma}

\begin{proof}
(i) is obvious, while (ii) follows immediately from
\cite[Lemma~3.2]{SeWe1}.
\end{proof}

\begin{proposition}\label{P:Letna}
Let $L$ be a lattice, let $\Sigma\subseteq\J(L)$. We consider the
following statements on $L$, $\Sigma$:
\begin{enumerate}
\item $L$ satisfies \Lt.

\item There are no elements $a$, $b$, $c$ of $\Sigma$ such that
$a\DD b\DD c$.
\end{enumerate}
Then \textup{(i)} implies \textup{(ii)}. Furthermore, if $\Sigma$ is
a join-seed of $L$, then \textup{(ii)} implies \textup{(i)}.
\end{proposition}

\begin{proof}
(i)$\Rightarrow$(ii) Suppose that there are $a$, $b$, $c\in\Sigma$
such that $a\DD b\DD c$. Let $b'$, $c'\in L$ such that both
inequalities $a\leq b\vee b'$ and $b\leq c\vee c'$ hold and are
minimal, respectively, in $b$ and in $c$. {}From the assumption that
$L$ satisfies \Lt\ it follows that
 \[
 a=(a\wedge b)\vee
 \Bigl(a\wedge\bigl((b\wedge c)\vee b'\bigr)\Bigr)\vee
 \Bigl(a\wedge\bigl((b\wedge c')\vee b'\bigr)\Bigr).
 \]
Since $a$ is \jirr\ and $a\nleq b$, there exists $x\in\set{c,c'}$
such that $a\leq(b\wedge x)\vee b'$. But $b\wedge x\leq b$, thus, by
the minimality statement on $b$, $b\leq x$, a contradiction.

(ii)$\Rightarrow$(i) under the additional assumption that $\Sigma$
is a join-seed of $L$. Let $a$, $b$, $b'$, $c$, $c'\in L$, denote by
$u$ (resp., $v$) the left hand side (resp., right hand side) of the
identity \Lt\ formed with these elements. It is clear that
$v\leq u$. Conversely, let $x\leq u$ in $\Sigma$, we prove that
$x\leq v$. If either $x\leq b\wedge(c\vee c')$ or $x\leq b'$ then
this is clear. Suppose that $x\nleq b\wedge(c\vee c'),b'$. Since
$x\leq\bigl(b\wedge(c\vee c')\bigr)\vee b'$ and $\Sigma$ is a
join-seed of $L$, there are $y\leq b\wedge(c\vee c')$ and $y'\leq b'$
in $\Sigma$ such that $x\leq y\vee y'$ is a minimal nontrivial
join-cover. If either $y\leq c$ or
$y\leq c'$ then either $x\leq a\wedge\bigl((b\wedge c)\vee b'\bigr)$
or $x\leq a\wedge\bigl((b\wedge c')\vee b'\bigr)$, in both cases
$x\leq v$. Suppose that $y\nleq c,c'$. Since $y\leq c\vee c'$ and
$\Sigma$ is a join-seed, there are $z\leq c$ and $z'\leq c'$ in
$\Sigma$ such that $y\leq z\vee z'$ is a minimal nontrivial
join-cover. Hence $x\DD y\DD z$, a contradiction. Therefore,
$x\leq v$. Since every element of $L$ is a join of elements of
$\Sigma$, $u\leq v$, whence $u=v$, which completes the proof that
$L$ satisfies \Lt.
\end{proof}

\begin{corollary}\label{C:Letna}
Let $(P,\utr)$ be a poset. Then $\Co(P)$ satisfies \Lt\ if{f}
$\lh P\leq 2$.
\end{corollary}

\begin{proof}
Put $\Sigma=\setm{\set{p}}{p\in P}$, the natural join-seed of
$\Co(P)$. Suppose first that $\lh P>2$, that is, $P$ contains a
four-element chain $o\tr a\tr b\tr c$. Then
$\set{a}\DD\set{b}\DD\set{c}$, thus, by Proposition~\ref{P:Letna},
$\Co(P)$ does not satisfy \Lt.

Conversely, suppose that $\Co(P)$ does not satisfy \Lt. By
Proposition~\ref{P:Letna}, there are $a$, $b$, $c\in P$ such that
$\set{a}\DD\set{b}\DD\set{c}$. Since $\set{a}\DD\set{b}$, there
exists $b'\in P$ such that either $b\tr a\tr b'$ or $b'\tr a\tr b$,
say, without loss of generality, $b'\tr a\tr b$. Since
$\set{b}\DD\set{c}$, there are $u$, $v\in P$ such that
$u\tr b\tr v$. Therefore, $b'\tr a\tr b\tr v$ is a four-element
chain in $P$.
\end{proof}

In order to proceed, it is convenient to recall the following result
from \cite{SeWe1}:

\begin{proposition}\label{P:UdBo}
Let $L$ be a complete, lower continuous, dually $2$-distributive
lattice that satisfies \Ud\ and \Bo. Then for every $p\in P$, there
are subsets $A$ and $B$ of $\rd{p}$ that satisfy the following
properties:
\begin{enumerate}
\item $\rd{p}=A\cup B$ and $A\cap B=\es$.

\item For all $x$, $y\in\rd{p}$, $p\leq x\vee y$ if{f} $(x,y)$ belongs
to $(A\times B)\cup(B\times A)$.
\end{enumerate}

Moreover, the set $\set{A,B}$ is uniquely determined by these
properties.
\end{proposition}

The set $\set{A,B}$ is called the \emph{Udav-Bond partition}
of $\rd{p}$ associated with~$p$.

We can now prove the following result:

\begin{theorem}\label{T:SUB2}
Let $L$ be a lattice. Then the following are equivalent:
\begin{enumerate}
\item $L$ belongs to $\SUB_2$.

\item $L$ satisfies the identities \Lt, \Ud, and \Bo.

\item There are a tree-like poset $\Gamma$ of length at most $2$
and a lattice embedding $\varphi\colon L\into\Co(\Gamma)$ that
preserves the existing bounds. Furthermore, the following additional
properties hold:
\begin{itemize}
\item[---] if $L$ is finite, then $\Gamma$ is finite;

\item[---] if $L$ is finite and subdirectly irreducible, then
$\varphi$ is atom-preserving.
\end{itemize}
\end{enumerate}
\end{theorem}

\begin{proof}
(i)$\Rightarrow$(ii)
It has been already proved in \cite{SeWe1} that every lattice in
$\SUB$ (thus \emph{a fortiori} in $\SUB_2$) satisfies the identities
\Ud\ and \Bo. Furthermore, it follows from Corollary~\ref{C:Letna}
that every lattice in $\SUB_2$ satisfies \Lt.

(ii)$\Rightarrow$(iii) Let $L$ be a lattice satisfying \Lt, \Ud, and
\Bo. We embed $L$ into the lattice $\hL=\Fil L$ of all filters of
$L$, partially ordered by reverse inclusion (see, e.g.,
G. Gr\"atzer~\cite{GLT}); if $L$ has no unit element, then we allow
the empty set in~$\hL$, otherwise we require filters to be nonempty.
This way, $\hL$ is a dually algebraic lattice, satisfies the same
identities as $L$, and the natural embedding $x\mapsto\upw x$
from~$L$ into $\hL$ preserves the existing bounds.

Hence we have reduced the problem to the case where $L$ is a dually
algebraic lattice. In particular, $L$ is complete, lower continuous,
and finitely spatial (it is even spatial), and $\Sigma=\J(L)$ is a
join-seed of $L$ (see Lemma~\ref{L:JoinSeeds}). Since $L$ satisfies
the identity \Lt\ and by Lemma~\ref{L:Lt22DD}, $L$ is dually
$2$-distributive. Hence, by Proposition~\ref{P:UdBo}, every
$p\in\J(L)$ has a unique Udav-Bond partition $\set{A_p,B_p}$.

Our poset $\Gamma$ is defined in a similar fashion as in
\cite[Section~7]{SeWe1}. The underlying set of $\Gamma$ is the set
of all nonempty finite sequences $\alpha=\seq{a_0,\dots, a_n}$ of
elements of $\J(L)$ such that $a_0$ is $\DD$-minimal in $\J(L)$
(this condition is added) and $a_i\DD a_{i+1}$, for all
$i\in\fso{n-1}$; as in \cite{SeWe1}, we call $n$ the \emph{length}
of $\alpha$ and we put $e(\alpha)=a_n$. Since $L$ satisfies \Lt\
and by Proposition~\ref{P:Letna}, the elements of $\Gamma$ are of
length either~$1$ or~$2$. Hence the partial ordering $\utr$ on
$\Gamma$ takes the following very simple form. The nontrivial
coverings in $\Gamma$ are those of the form
$\seq{p,a}\tr\seq{p}\tr\seq{p,b}$, where $p\in\J(L)$ and
$(a,b)\in A_p\times B_p$. Since the elements of length $1$ of
$\Gamma$ are either maximal or minimal, $\Gamma$ has indeed length at
most $2$. The proof that $\Gamma$ is tree-like proceeds
\emph{mutatis mutandis} as in \cite[Proposition~7.3]{SeWe1}.

As in \cite{SeWe1}, we define a map $\varphi$ from $L$ to the
powerset of $\Gamma$ by the rule
 \[
 \varphi(x)=\setm{\alpha\in\Gamma}{e(\alpha)\leq x},
 \qquad\text{for all }x\in L.
 \]
If $\seq{p,a}\tr\seq{p}\tr\seq{p,b}$ in $\Gamma$, then
$p\leq a\vee b$; hence, for $x\in L$, if both $\seq{p,a}$
and $\seq{p,b}$ belong to $\varphi(x)$, then $\seq{p}\in\varphi(x)$;
whence $\varphi(x)\in\Co(\Gamma)$.

It is clear that $\varphi$ is a \mh, and that it
preserves the existing bounds. Let $x$, $y\in L$ such that
$x\nleq y$. Since $L$ is finitely spatial, there exists
$a\in\nobreak\J(L)$ such that $a\leq x$ and $a\nleq y$. If $a$ is
$\DD$-minimal in $\J(L)$, then $\seq{a}$ belongs to
$\varphi(x)\setminus\varphi(y)$. If
$a$ is not $\DD$-minimal in $\J(L)$, then there exists $p\in\J(L)$
such that $p\DD a$. Since there are no $\DD$-chains with three
elements in $\J(L)$, $p$ is $\DD$-minimal, thus $\seq{p,a}$ belongs
to $\varphi(a)\setminus\varphi(b)$. Therefore, $\varphi$ is a
meet-embedding from $L$ into~$\Co(\Gamma)$.

We now prove that $\varphi$ is a \jh. It suffices to
prove that $\varphi(x\vee y)\subseteq\varphi(x)\vee\varphi(y)$, for
all $x$, $y\in L$. Let $\alpha\in\varphi(x\vee y)$, we prove that
$\alpha\in\varphi(x)\vee\varphi(y)$. This is obvious if
$\alpha\in\varphi(x)\cup\varphi(y)$, so suppose that
$\alpha\notin\varphi(x)\cup\varphi(y)$. Put $p=e(\alpha)$. So
$p\nleq x,y$ while $p\leq x\vee y$, thus there are $u\leq x$ and
$v\leq y$ in $\J(L)$ such that $p\leq u\vee v$ is a minimal
nontrivial join-cover. In particular, $p\DD u$ and $p\DD v$, thus
$\alpha=\seq{p}$ and both $\seq{p,u}$ and $\seq{p,v}$ belong to
$\Gamma$. It follows from $p\leq u\vee v$ that $(u,v)$ belongs to
$(A_p\times B_p)\cup(B_p\times A_p)$, thus either
$\seq{p,u}\tr\seq{p}\tr\seq{p,v}$ or
$\seq{p,v}\tr\seq{p}\tr\seq{p,u}$, in both cases
$\alpha\in\varphi(x)\vee\varphi(y)$. This completes the proof that
$\varphi$ is a lattice embedding.

Of course, if $L$ is finite, then $\Gamma$ is finite. Now suppose
that $L$ is finite and subdirectly irreducible. Since there are no
$\DD$-sequences of length three in $\J(L)$, there are \emph{a
fortiori} no $\DD$-cycles, thus, since $L$ is subdirectly
irreducible, $\J(L)$ has a unique $\DD$-minimal element $p$ (see
R. Freese, J. Je\v{z}ek, and J.\,B. Nation \cite[Chapter~3]{FJN}).
Hence, if $x$ is an atom of $L$, then $\varphi(x)$ is equal to
$\set{\seq{p}}$ if $x=p$ and to $\set{\seq{p,x}}$ otherwise, in both
cases, $\varphi(x)$ is an atom of $\Co(\Gamma)$.

Finally, (iii)$\Rightarrow$(i) is trivial.
\end{proof}

\begin{remark}
It follows from \cite[Example~8.1]{SeWe1} that there exists a
(non subdirectly irreducible) finite lattice $L$ without $\DD$-cycle
in $\SUB_2$ that cannot be embedded atom-preservingly into any
lattice of the form $\Co(P)$.
\end{remark}

\begin{proposition}\label{P:AtSUB2}
Let $L$ be a finite atomistic lattice without any $\DD$-cycle of the
form $a\DD b\DD a$. Then $L$ belongs to $\SUB$ if{f} $L$ belongs to
$\SUB_2$. In particular, $L$ has no $\DD$-cycle.
\end{proposition}

\begin{proof}
Suppose that $L$ belongs to $\SUB$. For $a$, $b$, $c\in\J(L)$
such that $a\DD b\DD c$, it follows from Lemma~\ref{L:JoinSeeds} that
there are elements $b'$ and $c'$ in $\J(L)$ such that both
inequalities $a\leq b\vee b'$ and
$b\leq c\vee c'$ hold and are minimal nontrivial join-covers. Since
$L$ satisfies \Stj, there exists $x\in\set{c,c'}$ such that
$b\leq a\vee x$ and $a\leq b'\vee x$. But $a\neq b$ and $b\neq x$
(because $a\DD b\DD x$), thus, since $a$, $b$, and $x$ are
atoms, the first inequality witnesses that $b\DD a$. Hence
$a\DD b\DD a$, a contradiction. It follows from
Proposition~\ref{P:Letna} that~$L$ satisfies \Lt, and then it
follows from Theorem~\ref{T:SUB2} that~$L$ belongs to $\SUB_2$, in
fact, there exists a finite poset $\Gamma$ of length at most $2$
such that $L$ embeds into $\Co(\Gamma)$. It follows from
Proposition~\ref{P:Letna} and Corollary~\ref{C:Letna} that
$\Co(\Gamma)$ has no $\DD$-cycle (a direct proof is also very easy),
thus neither has $L$.
\end{proof}

As the following example shows, Proposition~\ref{P:AtSUB2} does not
extend to the nonatomistic case.

\begin{example}\label{Ex:LBnonSUB2}
A finite subdirectly irreducible lattice without $\DD$-cycle that
belongs to $\SUB_3\setminus\SUB_2$.
\end{example}

\begin{proof}
Let $P=\set{\dot{a},\dot{a'},\dot{b},\dot{c},\dot{u},\dot{v}}$ be the
poset diagrammed on Figure~1.

\begin{figure}[htb]
\includegraphics{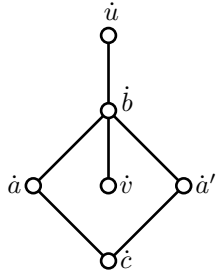}
\caption{A finite poset of length $3$}
\end{figure}

Let $L$ be the sublattice of $\Co(P)$
that consists of those subsets $X$ such that
 \begin{gather*}
 (\dot{a}\in X\Rightarrow\dot{a}'\in X)\text{ and }
 (\set{\dot{b},\dot{c}}\subseteq X\Rightarrow\dot{a}\in X)
 \text{ and }
 (\set{\dot{u},\dot{v}}\subseteq X\Rightarrow\dot{b}\in X)\\
 \text{ and }
 (\set{\dot{a}',\dot{u}}\subseteq X\Rightarrow\dot{b}\in X)
 \text{ and }
 (\set{\dot{u},\dot{c}}\subseteq X\Rightarrow\dot{a}\in X).
 \end{gather*}
Then $\J(L)=\set{a,a',b,c,u,v}$, where $a=\set{\dot{a},\dot{a}'}$,
$a'=\set{\dot{a}'}$, $b=\set{\dot{b}}$, $c=\set{\dot{c}}$,
$u=\set{\dot{u}}$, $v=\set{\dot{v}}$. Hence $L$ is the \jzs\ defined
by the generators $a$, $a'$, $b$, $c$, $u$, $v$, and the relations
 \[
 a'\leq a;\ a\leq b\vee c;\ b\leq u\vee v;\ b\leq a'\vee u;\
 a\leq u\vee c.
 \]
In particular, $L$ has no $\DD$-cycle and it is subdirectly
irreducible. Furthermore, $L$ is a sublattice of $\Co(P)$, hence it
belongs to $\SUB_3$. However, $L$ has the
three-element $\DD$-sequence $a\DD b\DD u$,
thus it does not belong to $\SUB_2$.
\end{proof}

\section{Local finiteness of $\SUB_2$}\label{S:LocFin2}

We begin with a few elementary observations on complete congruences
of lattices of the form $\Co(P)$. We recall that a congruence
$\theta$ of a complete lattice $L$ is \emph{complete}, if
$x\equiv y\pmod{\theta}$, for all $y\in Y$ implies
$x\equiv\bigvee Y\pmod{\theta}$ and
$x\equiv\bigwedge Y\pmod{\theta}$, for all $x\in L$ and all nonempty
$Y\subseteq L$. We say that $L$ is
\emph{completely subdirectly irreducible}, if it has a least nonzero
complete congruence.

\begin{definition}\label{D:Dclosed}
We say that a subset $U$ of a poset $(P,\utr)$ is
\emph{$\DD$-closed}, if $x\tr p\tr y$ and either $x\in U$ or $y\in U$
implies that $p\in U$, for all $x$, $y$, $p\in P$.
\end{definition}

Equivalently, $\set{p}\DD\set{x}$ (in $\Co(P)$) and
$x\in U$ implies that $p\in U$, for all $p$, $x\in P$.
Observe in particular that every $\DD$-closed subset of $P$ is
convex. We leave to the reader the straightforward proof of the
following lemma:

\begin{lemma}\label{L:P-U}
Let $P$ be a poset, let $U$ be a $\DD$-closed subset of $P$. Then
the binary relation $\theta_U$ on $\Co(P)$ defined by
 \[
 X\equiv Y\pmod{\theta_U}\Leftrightarrow X\cup U=Y\cup U,
 \quad\text{for all }X,\,Y\in\Co(P)
 \]
is a complete lattice congruence on $\Co(P)$, and one can define a
surjective homomorphism $h_U\colon\Co(P)\onto\Co(P\setminus U)$
with kernel $\theta_U$ by the rule $h_U(X)=X\setminus U$, for all
$X\in\Co(P)$. Furthermore, every complete lattice congruence
$\theta$ of $\Co(P)$ has the form $\theta_U$, with associated
$\DD$-closed set $U=\setm{p\in P}{\set{p}\equiv\es\pmod{\theta}}$.
\end{lemma}

We shall denote by $\CD(P)$ the lattice of all $\DD$-closed subsets
of a poset $P$ under inclusion. It follows from Lemma~\ref{L:P-U}
that \emph{$\CD(P)$ is isomorphic to the lattice of all complete
congruences of $\Co(P)$}.

\begin{lemma}\label{L:CDAlg}
The lattice $\CD(P)$ is algebraic, for every poset $P$.
\end{lemma}

\begin{proof}
Evidently, $\CD(P)$ is an algebraic subset of the powerset lattice
$\mathcal{P}(P)$ of $P$, that is, a complete
meet-subsemilattice closed under nonempty directed unions
(see~\cite{Gorb}). Since $\mathcal{P}(P)$ is algebraic, so
is $\CD(P)$.
\end{proof}

We observe that Lemma~\ref{L:CDAlg} cannot be extended to complete
congruences of arbitrary complete lattices: by G. Gr\"atzer and H.
Lakser~\cite{GrLa91}, every complete lattice~$L$ is isomorphic to
the lattice of complete congruences of some complete lattice $K$. By
G. Gr\"atzer and E.\,T. Schmidt~\cite{GrSc95}, $K$ can be taken
distributive.

\begin{corollary}\label{C:P-U}
For a poset $P$, the lattice $\Co(P)$ is completely subdirectly
irreducible if{f} there exists a least \pup{for the inclusion}
nonempty $\DD$-closed subset of $P$.
\end{corollary}

The analogue of Birkhoff's subdirect decomposition theorem
runs as follows:

\begin{lemma}\label{L:CplBirk}
Let $P$ be a poset. Then there exists a family $\famm{U_i}{i\in I}$
of $\DD$-closed subsets of $P$ such that the diagonal map from
$\Co(P)$ to $\prod_{i\in I}\Co(P\setminus U_i)$ is a lattice
embedding, and all the $\Co(P\setminus U_i)$ are completely
subdirectly irreducible.
\end{lemma}

\begin{proof}
Let $\setm{U_i}{i\in I}$ denote the set of all completely
meet-irreducible elements of $\CD(P)$. It follows from
Lemma~\ref{L:CDAlg} that $\CD(P)$ is dually spatial, that is,
every element of $\CD(P)$ is a meet of some of the $U_i$-s. By
applying this to the empty set, we obtain that the $U_i$-s have
empty intersection, which concludes the proof.
\end{proof}

\begin{notation}
For every positive integer $n$, we denote by $\bP_n$ the class of all
posets~$P$ of length at most $n$ such that $\Co(P)$ is completely
subdirectly irreducible (i.e., $P$ has a least nonempty
$\DD$-closed subset).
\end{notation}

For every pair $(I,J)$ of nonempty disjoint sets, set
$P_{I,J}=I\cup J\cup\set{p}$, where $p$ is some outside element, with
nontrivial coverings $x\tr p$ for $x\in I$ and
$p\tr y$ for $y\in J$.

\begin{lemma}\label{L:DescrP2}
The class $\bP_2$ consists of the one-element poset and all posets
of the form $P_{I,J}$, where $I$ and $J$ are nonempty disjoint sets.
\end{lemma}

\begin{proof}
It is straightforward to verify that the one-element poset and the
posets~$P_{I,J}$ all belong to $\bP_2$ (the monolith of
$\Co(P_{I,J})$ is the congruence $\Theta(\es,\set{p})$).
Conversely, let $P$ be a poset in $\bP_2$.
If $\lh P\leq 1$, then $\Co(P)$ is the powerset of $P$, thus it is
distributive. Furthermore, every subset of $P$
is $\DD$-closed, thus, since $P$ is completely subdirectly
irreducible,~$P$ is a singleton.

Suppose now that $P$ has length $2$. Thus there exists a
three-element chain $a\tr p\tr b$ in $P$. Since $P$ has length $2$,
$a$ is minimal, $b$ is maximal, and $\set{p}$ is $\DD$-closed. The
latter applies to every element of height $1$ instead of $p$, hence,
by assumption on $P$, $p$ is the only element of height $1$ of $P$.
Let $x$ be a minimal element of $P$. If $x\nutr p$, then $\set{x}$
is $\DD$-closed, thus $x=p$, a contradiction; whence $x\tr p$;
Similarly, $p\tr y$ for every maximal element $y$ of $P$. Therefore,
$P\cong P_{I,J}$, where $I$ (resp., $J$) is the set of all minimal
(resp., maximal) elements of $P$.
\end{proof}

\begin{notation}\label{Not:SUB2m}
For a positive integer $m$, let $\SUB_{2,m}$ denote the class of all
lattices that can be embedded into a product of lattices of the form
$\Co(P_{I,J})$, where $|I|+\nobreak|J|\leq m$.
\end{notation}

\begin{lemma}\label{L:truncCard}
Let $L$ be a finitely generated lattice, let $m\geq2$, let $a_0$,
\dots, $a_{m-1}$ be generators of $L$. Let $I$ and $J$ be disjoint
sets, let $f\colon L\to\Co(P_{I,J})$ be a lattice homomorphism. Then
there are finite sets $I'\subseteq I$ and $J'\subseteq J$ such that,
if $\pi\colon\Co(P_{I,J})\to\Co(P_{I',J'})$,
$X\mapsto X\cap P_{I',J'}$ is the canonical map, the following
assertions hold:
\begin{enumerate}
\item $|I'|+|J'|\leq 2^m-1$;

\item $\pi\circ f$ is a lattice homomorphism;

\item $\ker(f)=\ker(\pi\circ f)$.
\end{enumerate}
\end{lemma}

\begin{proof}
Let $\bD$ be the sublattice of the powerset lattice $\PP(I\cup J)$
generated by the subset $\setm{f(a_i)\setminus\set{p}}{i<m}$. We
observe that $\bD$ is a finite distributive lattice. Moreover, every
\jirr\ element of $\bD$ has the form $\bigwedge_{i\in X}f(a_i)$,
where $X$ is a proper subset of $\set{0,1,\dots,m-1}$, hence
$|\J(\bD)|\leq 2^m-2$.

\setcounter{claim}{0}
\begin{claim}\label{Cl:1}
The set $\bD^*=\Bigl(\bD\cap\bigl(\PP(I)\cup\PP(J)\bigr)\Bigr)\cup
\setm{X\cup\set{p}}{X\in\bD}$ is a sublattice of $\Co(P_{I,J})$, and
it contains the range of $f$.
\end{claim}

\begin{cproof}
It is easy to verify that $\bD^*$ is a sublattice of $\Co(P_{I,J})$.
It contains all elements of the form $f(a_i)$, thus it contains
the range of $f$.
\end{cproof}

For all $A\in\J(\bD)$, let $A^\dagger$ denote the largest element
$X$ of $\bD$ such that $A\not\subseteq X$. Observe that $A^\dagger$
is \mirr\ in $\bD$. For every $A\in\J(\bD)$, we pick
$k_A\in A\setminus A^\dagger$. Furthermore, if the zero $0_\bD$ of
$\bD$ is nonempty, we pick an element $l$ of $0_\bD$. We define
$K_0=\setm{k_A}{A\in\J(\bD)}$, and we put $K=K_0$ if $0_\bD=\es$,
$K=K_0\cup\set{l}$ otherwise. Observe that $K$ is a subset of
$I\cup J$ and $|K|\leq 2^m-1$. Finally, we put $I'=I\cap K$ and
$J'=J\cap K$, and we let $\pi\colon\Co(P_{I,J})\to\Co(P_{I',J'})$ be
the canonical map.

\begin{claim}\label{Cl:2}
The following assertions hold:
\begin{enumerate}
\item $X\not\subseteq Y$ implies that $X\cap K\not\subseteq Y\cap K$,
for all $X$, $Y\in\bD$.

\item $X\neq\es$ implies that $X\cap K\neq\es$, for all $X\in\bD$.
\end{enumerate}
\end{claim}

\begin{cproof}
(i) There exists $A\in\J(\bD)$ such that $A\subseteq X$ while
$A\not\subseteq Y$. Hence
$k_A\in A\setminus A^\dagger\subseteq X\setminus Y$.

(ii) If $0_\bD=\es$, then $X$ contains an atom $A$ of $\bD$; hence
$k_A\in A\subseteq X$. If $0_\bD\neq\es$, then
$l\in 0_\bD\subseteq X$.
\end{cproof}

Now we can prove that $\pi\circ f$ is a lattice homomorphism. It is
clearly a \mh. To prove that it is a \jh, it suffices to prove the
containment
 \begin{equation}\label{Eq:pifhom}
 (f(x)\vee f(y))\cap P_{I',J'}\subseteq
 (f(x)\cap P_{I',J'})\vee(f(y)\cap P_{I',J'}),
 \end{equation}
for all $x$, $y\in L$. Suppose otherwise. Since $p$ is the only
element of $P_{I,J}$ that is neither maximal nor minimal, it
belongs to the left hand side of \eqref{Eq:pifhom} but not
to its right hand side.
In particular, $p\notin f(x)\cup f(y)$, whence, say,
$f(x)\subseteq I$ and $f(y)\subseteq J$. By Claim~\ref{Cl:1},
$f(x)$, $f(y)\in\bD^*$, thus $f(x)$, $f(y)\in\bD$. Furthermore,
$p\in f(x)\vee f(y)$ with $f(x)\subseteq I$ and $f(y)\subseteq J$,
whence $f(x)$, $f(y)$ are nonempty. By Claim~2(ii), both $f(x)$ and
$f(y)$ meet $K$, whence $p\in(f(x)\cap I')\vee(f(y)\cap J')$, a
contradiction. Therefore, $\pi\circ f$ is indeed a lattice
homomorphism.

In order to conclude the proof of Lemma~\ref{L:truncCard}, it
suffices to prove that $\ker(\pi\circ f)$ is contained in $\ker(f)$.
So let $x$, $y\in L$ such that $f(x)\not\subseteq f(y)$. By
Claim~\ref{Cl:1}, both $f(x)$ and $f(y)$ belong to $\bD^*$. If
$f(x)\setminus\set{p}\subseteq f(y)$, then $p\in f(x)$, hence
 \[
 p\in\bigl(f(x)\cap P_{I',J'}\bigr)\setminus
 \bigl(f(y)\cap P_{I',J'}\bigr)=
 (\pi\circ f(x))\setminus(\pi\circ f(y)).
 \]
If $f(x)\setminus\set{p}\not\subseteq f(y)$, then, by Claim~2(i),
there exists $k\in K$ with
$k\in(f(x)\setminus\set{p})\setminus(f(y)\setminus\set{p})$, whence
$k\in(\pi\circ f(x))\setminus(\pi\circ f(y))$. In both cases,
$\pi\circ f(x)\not\subseteq\pi\circ f(y)$.
\end{proof}

We can now prove the main result of this section:

\begin{theorem}\label{T:SUB2lf}
Let $m\geq2$ be an integer. Then every $m$-generated member of
$\SUB_2$ belongs to $\SUB_{2,2^m-1}$. In particular, the variety
$\SUB_2$ is locally finite.
\end{theorem}

\begin{proof}
Let $L$ be a $m$-generated member of $\SUB_2$.
By Lemma~\ref{L:CplBirk}, there exists a family
$\famm{(I_l,J_l)}{l\in\Omega}$ of pairs of nonempty disjoint sets,
together with an embedding
$f\colon L\into\prod_{l\in\Omega}\Co(P_{I_l,J_l})$. For all
$l\in\Omega$, denote by $f_l\colon L\to\Co(P_{I_l,J_l})$ the $l$-th
component of $f$.
By Lemma~\ref{L:truncCard}, there are finite subsets
$I'_l\subseteq I_l$ and $J'_l\subseteq J_l$ such that
$|I'_l|+|J'_l|\leq 2^m-1$, $\pi_l\circ f_l$ is a lattice
homomorphism, and $\ker(f_l)=\ker(\pi_l\circ f_l)$, where
$\pi_l\colon\Co(P_{I_l,J_l})\to\Co(P_{I'_l,J'_l})$ is the canonical
map. Therefore, the map
 \[
 g\colon L\to\prod_{l\in\Omega}\Co(P_{I'_l,J'_l}),\
 x\mapsto\famm{\pi_l\circ f_l(x)}{l\in\Omega}
 \]
is a lattice embedding of $L$ into a member of $\SUB_{2,2^m-1}$.
\end{proof}

The above argument gives a very rough upper bound for the
cardinality of the free lattice $F_m$ in $\SUB_2$ on $m$ generators,
namely, $e(m)^{e(m)^m}$, where $e(m)=2^{2^m}+2^{2^{m+1}-2}-1$. Indeed,
by Theorem~\ref{T:SUB2lf}, $F_m$ embeds into $A^{A^m}$, where
$A=P_{2^m-1,2^m-1}$, and $\vert A\vert=e(m)$.

\section{The identities \Ht{n}}\label{S:Hn}

\begin{definition}\label{D:Pols}
For a positive integer $n$, we define inductively lattice
polynomials $U_{i,n}$ (for $0\leq i\leq n$), $V_{i,j,n}$ (for
$0\leq j\leq i\leq n-1$), $W_{i,j,n}$ (for $0\leq j\leq i\leq n-2$),
with variables $x_0$, \dots, $x_n$, $x'_1$, \dots, $x'_n$, as
follows:
 \begin{align*}
 U_{n,n}&=x_n;\\
 U_{i,n}&=x_i\wedge(U_{i+1,n}\vee x'_{i+1})
 &&\text{for }0\leq i\leq n-1;\\
 V_{i,i,n}&=(x_i\wedge U_{i+1,n})\vee(x_i\wedge x'_{i+1})
 &&\text{for }0\leq i\leq n-1;\\
 V_{i,j,n}&=x_j\wedge(V_{i,j+1,n}\vee x'_{j+1})
 &&\text{for }0\leq j<i\leq n-1;\\
 W_{i,i,n}&=x_i\wedge(x'_{i+1}\vee x'_{i+2})\wedge
 \bigl((U_{i+1,n}\wedge(x_i\vee x'_{i+2}))\vee x'_{i+1}\bigr)
 &&\text{for }0\leq i\leq n-2;\\
 W_{i,j,n}&=x_j\wedge(W_{i,j+1,n}\vee x'_{j+1})
 &&\text{for }0\leq j<i\leq n-2.
 \end{align*}
Furthermore, we put
 \begin{align*}
 U_n&=U_{0,n},\\
 V_{i,n}&=V_{i,0,n}&\text{for }0\leq i\leq n-1;\\
 W_{i,n}&=W_{i,0,n}&\text{for }0\leq i\leq n-2.
 \end{align*}
\end{definition}

\begin{lemma}\label{L:TrivIneq}
Let $n$ be a positive integer. The following inequalities hold in
every lattice:
\begin{enumerate}
\item $V_{i,j,n}\leq U_{j,n}$ for $0\leq j\leq i\leq n-1$;
\item $W_{i,j,n}\leq U_{j,n}$ for $0\leq j\leq i\leq n-2$;
\item $V_{i,n}\leq U_n$ for $0\leq i\leq n-1$;
\item $W_{i,n}\leq U_n$ for $0\leq i\leq n-2$.
\end{enumerate}
\end{lemma}

\begin{proof}
Items (i) and (ii) are easily established by downward induction on
$j$. Items (iii) and (iv) follow immediately.
\end{proof}

As in the following lemma, we shall often use the convenient notation
 \[
 \vec a=\seq{a_0,a_1,\ldots,a_n},\qquad
 \vec a'=\seq{a'_1,\ldots,a'_n}.
 \]

\begin{lemma}\label{L:Id2JIT1}
Let $n$ be a positive integer, let $L$ be a lattice, let $a_0$,
\dots, $a_n\in\J(L)$ and $a'_1$, \dots, $a'_n\in L$ such that
$a_i\leq a_{i+1}\vee a'_{i+1}$ is a nontrivial join-cover, for all
$i\in\fso{n-1}$, minimal in $a_{i+1}$ for $i\leq n-2$.
If the equality
 \begin{equation}\label{Eq:a0bigmess}
 a_0=\bigvee_{0\leq i\leq n-1}V_{i,n}(\vec a,\vec a')\vee
\bigvee_{0\leq i\leq n-2}W_{i,n}(\vec a,\vec a')
 \end{equation}
holds, then there exists $i\in\fso{n-2}$ such that
$a_i\leq a'_{i+1}\vee a'_{i+2}$ and $a_{i+1}\leq a_i\vee a'_{i+2}$.
\end{lemma}

\begin{note}
Of course, the meaning of the right hand side of the equation
\eqref{Eq:a0bigmess} for $n=1$ is simply $V_{0,1}(\vec a,\vec a')$.
\end{note}

\begin{proof}
We first observe that the assumptions imply the following:
 \begin{equation}\label{Eq:U=a}
 U_{i,n}(\vec a,\vec a')=a_i,\text{ for all }i\in\fso{n}.
 \end{equation}
Now we put $c_{i,j}=V_{i,j,n}(\vec a,\vec a')$ and $c_i=c_{i,0}$ for
$0\leq j\leq i\leq n-1$, and $d_{i,j}=W_{i,j,n}(\vec a,\vec a')$ and
$d_i=d_{i,0}$ for $0\leq j\leq i\leq n-2$. We deduce from the
assumption that one of the two following cases occurs:

\begin{itemize}
\item[\textbf{Case 1.}] $a_0=c_i$ for some $i\in\fso{n-1}$. This can
also be written $c_{i,0}=\nobreak a_0$. Suppose that $c_{i,j}=a_j$,
for $0\leq j<i$. So $a_j\leq c_{i,j+1}\vee a'_{j+1}$ with
$c_{i,j+1}\leq a_{j+1}$, thus, by the minimality assumption on
$a_{j+1}$, we obtain that $c_{i,j+1}=\nobreak a_{j+1}$. Hence
$c_{i,j}=a_j$, for all $j\in\fso{i}$, in particular,
by~\eqref{Eq:U=a},
 \[
 a_i=c_{i,i}=(a_i\wedge a_{i+1})\vee(a_i\wedge a'_{i+1}),
 \]
whence, by the \jirry\ of $a_i$, either $a_i\leq a_{i+1}$ or
$a_i\leq a'_{i+1}$, which contradicts the assumption.
Thus, Case~1 cannot occur.

\item[\textbf{Case 2.}] $a_0=d_i$ for some $i\in\fso{n-2}$ (thus
$n\geq2$). As in Case~1, $d_{i,j}=a_j$, for all $j\in\fso{i}$,
whence, for $j=i$ and by \eqref{Eq:U=a},
 \[
 a_i\leq(a'_{i+1}\vee a'_{i+2})\wedge
 \bigl((a_{i+1}\wedge(a_i\vee a'_{i+2}))\vee a'_{i+1}\bigr)
 \]
Set $x=a_{i+1}\wedge(a_i\vee a'_{i+2})$, so
$x\leq a_{i+1}$. Observe that $a_i\leq a'_{i+1}\vee a'_{i+2}$
and $a_i\leq x\vee a'_{i+1}$,
whence, by the minimality assumption on $a_{i+1}$, we obtain that
$x=a_{i+1}$, that is, $a_{i+1}\leq a_i\vee a'_{i+2}$.
\end{itemize}
This concludes the proof.
\end{proof}

\begin{lemma}\label{L:JIT2Id1}
Let $L$ be a lattice satisfying the Stirlitz identity \St, let
$\Sigma$ be a join-seed of $L$, let $x\in\Sigma$, let
$n$ be a positive integer, and let $a_0$, \dots, $a_n$, $a'_1$,
\dots, $a'_n\in L$. If $x\leq U_n(\vec a,\vec a')$, then one of the
following three cases occurs:
\begin{enumerate}
\item there exists $i\in\fso{n-1}$ such that
$x\leq V_{i,n}(\vec a,\vec a')$;

\item there exists $i\in\fso{n-2}$ such that
$x\leq W_{i,n}(\vec a,\vec a')$;

\item there are elements $x_i\leq U_{i,n}(\vec a,\vec a')$
\pup{$0\leq i\leq n$} and $x'_i\leq a'_i$ \pup{$1\leq i\leq n$}
of $\Sigma$ such that the pair
$(\famm{x_i}{0\leq i\leq n},\famm{x'_i}{1\leq i\leq n})$
is a Stirlitz track.
\end{enumerate}
\end{lemma}

\begin{proof}
We put $a^*_i=U_{i,n}(\vec a,\vec a')$ for $0\leq i\leq n$,
$c_{i,j}=V_{i,j,n}(\vec a,\vec a')$ for $0\leq j\leq i\leq n-1$ and
$d_{i,j}=W_{i,j,n}(\vec a,\vec a')$ for $0\leq j\leq i\leq n-2$,
then $c_i=c_{i,0}$ for $0\leq i\leq n-1$ and $d_i=d_{i,0}$ for
$0\leq i\leq n-2$. We observe that
$x\leq U_{0,n}(\vec a,\vec a')=a^*_0$.

Suppose that $x\nleq c_i$, for all $i\in\fso{n-1}$. Put $x_0=x$.
Suppose we have constructed $x_j\leq a^*_j$ in $\Sigma$, with
$0\leq j<n$, such that $x_j\nleq c_{i,j}$, for all
$i\in\set{j,\dots,n-1}$. If either $x_j\leq a^*_{j+1}$ or
$x_j\leq a'_{j+1}$, then, since $x_j\leq a_j$, we obtain that
$x_j\leq c_{j,j}$, a contradiction; whence
$x_j\nleq a^*_{j+1},a'_{j+1}$. On the other hand,
$x_j\leq a^*_j\leq a^*_{j+1}\vee a'_{j+1}$, thus, since
$x_j\in\Sigma$ and $\Sigma$ is a join-seed of $L$, there
are $x_{j+1}\leq a^*_{j+1}$ and $x'_{j+1}\leq a'_{j+1}$ in $\Sigma$
such that $x_j\leq x_{j+1}\vee x'_{j+1}$ is a minimal nontrivial
join-cover. Suppose that $x_{j+1}\leq c_{i,j+1}$ for some
$i\in\set{j+1,\dots,n-1}$. Then
 \[
 x_j\leq a_j\wedge(x_{j+1}\vee x'_{j+1})\leq
 a_j\wedge(c_{i,j+1}\vee a'_{j+1})=c_{i,j},
 \]
a contradiction. Hence $x_{j+1}\nleq c_{i,j+1}$, for all
$i\in\set{j+1,\dots,n-1}$, which completes the induction step.

Therefore, we have constructed elements $x_0\leq a^*_0$, \dots,
$x_n\leq a^*_n$, $x'_1\leq a'_1$, \dots, $x'_n\leq a'_n$ of
$\Sigma$ such that $x_0=x$ and $x_i\leq x_{i+1}\vee x'_{i+1}$ is a
minimal nontrivial join-cover, for all $i\in\fso{n-1}$.
Suppose that $(\famm{x_i}{0\leq i\leq n},\famm{x'_i}{1\leq i\leq n})$
is not a Stirlitz track. Then, since all the $x_i$-s and the
$x'_i$-s are \jirr\ and~$L$ satisfies the axiom \Stj\ (see
\cite[Proposition~4.4]{SeWe1}), there exists
$i\in\fso{n-2}$ such that
 \begin{equation}\label{Eq:nonSti}
 x_{i+1}\leq x_i\vee x'_{i+2}\text{ and }
 x_i\leq x'_{i+1}\vee x'_{i+2}.
 \end{equation}
It follows from this that $x_{i+1}\leq a^*_{i+1}\wedge(a_i\vee
a'_{i+2})$, whence
 \[
 x_i\leq a_i\wedge(a'_{i+1}\vee a'_{i+2})\wedge
 \bigl((a^*_{i+1}\wedge(a_i\vee a'_{i+2}))\vee a'_{i+1}\bigr)
 =d_{i,i}.
 \]
For $0\leq j<i$, suppose we have proved that $x_{j+1}\leq d_{i,j+1}$.
Since $x_j\leq x_{j+1}\vee x'_{j+1}$, we obtain that
$x_j\leq a_j\wedge(d_{i,j+1}\vee a'_{j+1})=d_{i,j}$. Hence we have
proved that $x_j\leq d_{i,j}$, for all $j\in\fso{i}$. In particular,
$x=x_0\leq d_{i,0}=d_i=W_{i,n}(\vec a,\vec a')$, which concludes the
proof.
\end{proof}

For a positive integer $n$, let \Ht{n} be the following lattice
identity:
 \[
 U_n=\bigvee_{0\leq i\leq n-1}V_{i,n}\vee
 \bigvee_{0\leq i\leq n-2}W_{i,n}.
 \]
It is not hard to verify directly that \Ht{1} is equivalent to
distributivity.

\begin{proposition}\label{P:MainHn}
Let $n$ be a positive integer, let $L$ be a lattice satisfying \St\
and~\Ud, let $\Sigma$ be a subset of $\J(L)$. We consider the
following statements on $L$, $\Sigma$:
\begin{enumerate}
\item $L$ satisfies \Ht{n}.

\item For all elements $a_0$, \dots, $a_n$, $a'_1$, \dots, $a'_n$ of
$\Sigma$, if $a_i\leq a_{i+1}\vee a'_{i+1}$ is a nontrivial
join-cover, for all $i\in\fso{n-1}$, minimal in $a_{i+1}$ for
$i\neq n-1$, then there exists $i\in\fso{n-2}$ such that
$a_i\leq a'_{i+1}\vee a'_{i+2}$ and $a_{i+1}\leq a_i\vee a'_{i+2}$.

\item There is no Stirlitz track of length $n$ with entries in
$\Sigma$.
\end{enumerate}
Then \textup{(i)} implies \textup{(ii)} implies \textup{(iii)}.
Furthermore, if $\Sigma$ is a join-seed of $L$, then \textup{(iii)}
implies \textup{(i)}.
\end{proposition}

\begin{proof}
(i)$\Rightarrow$(ii) Let $a_0$, \dots, $a_n$, $a'_1$,
\dots, $a'_n\in\Sigma$ satisfy the assumption of (ii). Observe that
$U_{i,n}(\vec a,\vec a')=a_i$ for $0\leq i\leq n$, in particular,
$U_n(\vec a,\vec a')=a_0$. {}From the assumption that $L$ satisfies
\Ht{n} it follows that
 \[
 a_0=\bigvee_{0\leq i\leq n-1}V_{i,n}(\vec a,\vec a')\vee
 \bigvee_{0\leq i\leq n-2}W_{i,n}(\vec a,\vec a').
 \]
The conclusion of (ii) follows from Lemma~\ref{L:Id2JIT1}.

(ii)$\Rightarrow$(iii)
Let $\sigma=(\famm{a_i}{0\leq i\leq n},\famm{a'_i}{1\leq i\leq n})$
be a Stirlitz track of $L$ with entries in $\Sigma$. {}From (ii)
it follows that there exists $i\in\fso{n-2}$ such that 
$a_i\leq a'_{i+1}\vee a'_{i+2}$ and $a_{i+1}\leq a_i\vee a'_{i+2}$,
whence $a_{i+1}\leq a'_{i+1}\vee a'_{i+2}$. Since $\sigma$ is a
Stirlitz track, the inequality $a_{i+1}\leq a'_{i+1}\vee a_{i+2}$
also holds, whence, since $a_{i+1}\leq a_{i+2}\vee a'_{i+2}$ and by
\Udj, either $a_{i+1}\leq a'_{i+1}$ or $a_{i+1}\leq a_{i+2}$ or
$a_{i+1}\leq a'_{i+2}$, a contradiction.

(iii)$\Rightarrow$(i) under the additional assumption that $\Sigma$
is a join-seed of $L$. Let $a_0$, \dots, $a_n$, $a'_1$, \dots,
$a'_n\in L$, define $c$, $d\in L$ by
 \[
 c=U_n(\vec a,\vec a'),\qquad
 d=\bigvee_{0\leq i\leq n-1}V_{i,n}(\vec a,\vec a')\vee
 \bigvee_{0\leq i\leq n-2}W_{i,n}(\vec a,\vec a').
 \]
It follows from Lemma~\ref{L:TrivIneq} that $d\leq c$. Conversely,
let $x\in\Sigma$ such that $x\leq c$, we prove that $x\leq d$.
Otherwise, $x\nleq V_{i,n}(\vec a,\vec a')$, for all $i\in\fso{n-1}$
and $x\nleq W_{i,n}(\vec a,\vec a')$, for all $i\in\fso{n-2}$, thus,
by Lemma~\ref{L:JIT2Id1}, there are elements $x_0=x$, $x_1$, \dots,
$x_n$, $x'_1$, \dots, $x'_n$ of $\Sigma$ such that the pair
 \[
 (\famm{x_i}{0\leq i\leq n},\famm{x'_i}{1\leq i\leq n})
 \]
is a Stirlitz track of $L$, a contradiction. Since every element of
$L$ is a join of elements of $\Sigma$, it follows that $c\leq d$.
Therefore, $c=d$, so $L$ satisfies \Ht{n}.
\end{proof}

\begin{corollary}\label{C:MainHnP}
Let $(P,\utr)$ be a poset, let $n$ be a positive integer. Then
$\Co(P)$ satisfies \Ht{n} if{f} $\lh P\leq n$.
\end{corollary}

\begin{proof}
It follows from \cite[Section~4]{SeWe1}
that $\Co(P)$ satisfies \St\ and \Ud. Furthermore,
$\Sigma=\setm{\set{p}}{p\in P}$ is a join-seed of $\Co(P)$.

Suppose first that $\lh P\geq n+1$, that is, $P$ contains a
$n+2$-element chain, say, $y\tr x_0\tr\cdots\tr x_n$. Then the pair
 \[
 (\famm{\set{x_i}}{0\leq i\leq n},\famm{\set{y}}{1\leq i\leq n})
 \]
is a Stirlitz track of length $n$ in $\Co(P)$, thus, by
Proposition~\ref{P:MainHn}, $\Co(P)$ does not satisfy \Ht{n}.

Conversely, suppose that $P$ does not contain any $n+2$-element
chain. By Proposition~\ref{P:MainHn}, in order to prove that
$\Co(P)$ satisfies \Ht{n}, it suffices to prove that $\Co(P)$ has no
Stirlitz track of length $n$ with entries in $\Sigma$. Suppose that
there exists such a Stirlitz track, say,
 \[
 (\famm{\set{x_i}}{0\leq i\leq n},\famm{\set{x'_i}}{1\leq i\leq n}).
 \]
Since $\set{x_0}\leq\set{x_1}\vee\set{x'_1}$ is a nontrivial
join-cover, either $x_1\tr x_0\tr x'_1$ or $x'_1\tr x_0\tr x_1$,
say, $x'_1\tr x_0\tr x_1$. Similarly, for all $i\in\fso{n-1}$,
either $x_{i+1}\tr x_i\tr x'_{i+1}$ or $x'_{i+1}\tr x_i\tr x_{i+1}$.
Suppose that the first possibility occurs, and take~$i$ minimum
such. Thus $i>0$ and $x'_i\tr x_{i-1}\tr x_i\tr x'_{i+1}$ and
$x_{i+1}\tr x_i$ while $\set{x_i}\leq\set{x'_i}\vee\set{x_{i+1}}$, a
contradiction. Thus $x'_{i+1}\tr x_i\tr x_{i+1}$. It follows that
 \[
 x'_1\tr x_0\tr\cdots\tr x_n
 \]
is a $n+2$-element chain in $P$, a contradiction.
\end{proof}

\section{The identities \Ht{m,n}}\label{S:Hmn}

\begin{definition}\label{D:bi-Stirlitz}
For positive integers $m$ and $n$ and a lattice $L$, a
\emph{bi-Stirlitz track of index $(m,n)$} is a pair $(\sigma,\tau)$,
where
 \begin{align*}
 \sigma&=(\famm{a_i}{0\leq i\leq m},\famm{a'_i}{1\leq i\leq m}),\\
 \tau&=(\famm{b_j}{0\leq j\leq n},\famm{b'_j}{1\leq j\leq n})
 \end{align*}
are Stirlitz tracks with the same base $a_0=b_0\leq a_1\vee b_1$.
\end{definition}

For positive integers $m$ and $n$, we define the identity \Ht{m,n},
with variable symbols $t$, $x_i$, $x'_i$ ($1\leq i\leq m$), $y_j$,
$y'_j$ ($1\leq j\leq n$) as follows, where we put $x_0=y_0=t$:
 \begin{align*}
 U_m(\vec x,\vec x')\wedge U_n(\vec y,\vec y')=&
 \bigvee_{0\leq i\leq m-1}\bigl(
 V_{i,m}(\vec x,\vec x')\wedge U_n(\vec y,\vec y')\bigr)\\
 &\vee\bigvee_{0\leq i\leq m-2}\bigl(
 W_{i,m}(\vec x,\vec x')\wedge U_n(\vec y,\vec y')\bigr)\\
 &\vee\bigvee_{0\leq j\leq n-1}\bigl(
 U_m(\vec x,\vec x')\wedge V_{j,n}(\vec y,\vec y')\bigr)\\
 &\vee\bigvee_{0\leq j\leq n-2}\bigl(
 U_m(\vec x,\vec x')\wedge W_{j,n}(\vec y,\vec y')\bigr)\\
 &\vee\bigl(U_m(\vec x,\vec x')\wedge U_n(\vec y,\vec y')
 \wedge(x_1\vee y'_1)\wedge(x'_1\vee y_1)
 \bigr).
 \end{align*}

The analogue of Proposition~\ref{P:MainHn} for the identity \Ht{m,n}
is the following:

\begin{proposition}\label{P:MainHmn}
Let $m$ and $n$ be positive integers, let $L$ be a lattice satisfying
\St, \Ud, and \Bo, let $\Sigma$ be a subset of $\J(L)$. We consider
the following statements on $L$, $\Sigma$:
\begin{enumerate}
\item $L$ satisfies \Ht{m,n}.

\item For all elements $a_0$, \dots, $a_m$, $a'_1$, \dots, $a'_m$,
$b_0$, \dots, $b_n$, $b'_1$, \dots, $b'_n$ of
$\Sigma$ with $a_0=b_0$, if $a_i\leq a_{i+1}\vee a'_{i+1}$ is a
nontrivial join-cover, for all $i\in\fso{m-1}$, minimal in $a_{i+1}$
for $i\neq m-1$ and if $b_j\leq b_{j+1}\vee b'_{j+1}$ is a
nontrivial join-cover, for all $j\in\fso{n-1}$, minimal in $b_{j+1}$
for $j\neq n-1$, then one of the following occurs:
 \begin{enumerate}
 \item there exists $i\in\fso{m-2}$ such that
 $a_i\leq a'_{i+1}\vee a'_{i+2}$ and\linebreak
$a_{i+1}\leq a_i\vee a'_{i+2}$;

 \item there exists $j\in\fso{n-2}$ such that
 $b_j\leq b'_{j+1}\vee b'_{j+2}$ and\linebreak
$b_{j+1}\leq b_j\vee b'_{j+2}$;

 \item $a_0\leq(a_1\vee b'_1)\wedge(a'_1\vee b_1)$.
 \end{enumerate}
\item There is no bi-Stirlitz track of index $(m,n)$ with entries in
$\Sigma$.
\end{enumerate}

Then \textup{(i)} implies \textup{(ii)} implies \textup{(iii)}.
Furthermore, if $\Sigma$ is a join-seed of $L$, then \textup{(iii)}
implies \textup{(i)}.
\end{proposition}

\begin{proof}
(i)$\Rightarrow$(ii) Let $a_0$, \dots, $a_m$, $a'_1$,
\dots, $a'_m$, $b_0$, \dots, $b_n$, $b'_1$, \dots, $b'_n\in\Sigma$
satisfy the assumption of (ii). Observe that
$U_{m,i}(\vec a,\vec a')=a_i$ for $0\leq i\leq m$ and
$U_{n,j}(\vec b,\vec b')=b_j$ for $0\leq j\leq n$. Put $p=a_0=b_0$.
{}From the assumption that $L$ satisfies \Ht{m,n} it follows that
\begin{equation}\label{Eq:p=bigmess}
 \begin{aligned}
 p=&
 \bigvee_{0\leq i\leq m-1}\bigl(
 V_{i,m}(\vec a,\vec a')\wedge U_n(\vec b,\vec b')\bigr)
 \vee\bigvee_{0\leq i\leq m-2}\bigl(
 W_{i,m}(\vec a,\vec a')\wedge U_n(\vec b,\vec b')\bigr)\\
 &\vee\bigvee_{0\leq j\leq n-1}\bigl(
 U_m(\vec a,\vec a')\wedge V_{j,n}(\vec b,\vec b')\bigr)
 \vee\bigvee_{0\leq j\leq n-2}\bigl(
 U_m(\vec a,\vec a')\wedge W_{j,n}(\vec b,\vec b')\bigr)\\
 &\vee\bigl(U_m(\vec a,\vec a')\wedge U_n(\vec b,\vec b')
 \wedge(a_1\vee b'_1)\wedge(a'_1\vee b_1)
 \bigr).
 \end{aligned}
\end{equation}
Since $p$ is \jirr, three cases can occur:
\begin{itemize}
\item[\textbf{Case 1.}]
$p=\bigvee\limits_{0\leq i\leq m-1}\bigl(
V_{i,m}(\vec a,\vec a')\wedge U_n(\vec b,\vec b')\bigr)
\vee\bigvee\limits_{0\leq i\leq m-2}\bigl(
W_{i,m}(\vec a,\vec a')\wedge U_n(\vec b,\vec b')\bigr)$.

{}From Lemma~\ref{L:TrivIneq} it follows that the equality
 \[
 p=\bigvee_{0\leq i\leq m-1}V_{i,m}(\vec a,\vec a')
 \vee\bigvee_{0\leq i\leq m-2}W_{i,m}(\vec a,\vec a')
 \]
also holds. By Lemma~\ref{L:Id2JIT1}, there exists $i\in\fso{m-2}$
such that $a_i\leq a'_{i+1}\vee a'_{i+2}$ and
$a_{i+1}\leq a_i\vee a'_{i+2}$.

\item[\textbf{Case 2.}]
$p=\bigvee\limits_{0\leq j\leq n-1}\bigl(
U_m(\vec a,\vec a')\wedge V_{j,n}(\vec b,\vec b')\bigr)
\vee\bigvee\limits_{0\leq j\leq n-2}\bigl(
U_m(\vec a,\vec a')\wedge W_{j,n}(\vec b,\vec b')\bigr)$.

As in Case~1, we obtain $j\in\fso{n-2}$
such that $b_j\leq b'_{j+1}\vee b'_{j+2}$ and
$b_{j+1}\leq b_j\vee b'_{j+2}$.

\item[\textbf{Case 3.}]
$p\leq(a_1\vee b'_1)\wedge(a'_1\vee b_1)$.
\end{itemize}

In all three cases above, the conclusion of (ii) holds.

(ii)$\Rightarrow$(iii) Let $(\sigma,\tau)$ be a bi-Stirlitz track as
in Definition~\ref{D:bi-Stirlitz}. Put $p=a_0=b_0$.
It follows from the assumption
(ii) that either there exists $i\in\fso{m-2}$ such that
$a_i\leq a'_{i+1}\vee a'_{i+2}$ and
$a_{i+1}\leq a_i\vee a'_{i+2}$, or there exists $j\in\fso{n-2}$
such that $b_j\leq b'_{j+1}\vee b'_{j+2}$ and
$b_{j+1}\leq b_j\vee b'_{j+2}$, or
$p\leq(a_1\vee b'_1)\wedge(a'_1\vee b_1)$. In the first case,
$a_{i+1}\leq a'_{i+1}\vee a'_{i+2}$, but $\sigma$ is a Stirlitz
track, thus also $a_{i+1}\leq a'_{i+1}\vee a_{i+2}$, a contradiction
since $a_{i+1}\leq a_{i+2}\vee a'_{i+2}$ and by \Udj. The second
case leads to a similar contradiction. In the third case,
$p\leq a_1\vee b'_1$, a contradiction by \Udj\ since
$p\leq a_1\vee b_1$ and $p\leq a_1\vee a'_1$.

(iii)$\Rightarrow$(i) under the additional assumption that $\Sigma$
is a join-seed of $L$. Let $a_0=b_0$, $a_1$, \dots, $a_m$, $a'_1$,
\dots, $a'_m$, $b_1$, \dots, $b_n$, $b'_1$, \dots, $b'_n\in L$,
put $c=U_m(\vec a,\vec a')\wedge U_n(\vec b,\vec b')$ and define
$d\in L$ as the right hand side
of \eqref{Eq:p=bigmess}. Further, put
$a^*_i=U_{i,m}(\vec a,\vec a')$ for $0\leq i\leq m$ and 
$b^*_j=U_{j,n}(\vec b,\vec b')$ for $0\leq j\leq n$.
It follows from Lemma~\ref{L:TrivIneq} that
$d\leq c$. Conversely, let $z\in\Sigma$ such that $z\leq c$, we
prove that $z\leq d$. Otherwise, $z\nleq V_{i,m}(\vec a,\vec a')$,
for all $i\in\fso{m-1}$, and $z\nleq W_{i,m}(\vec a,\vec a')$, for
all $i\in\fso{m-2}$, and $z\nleq V_{j,n}(\vec b,\vec b')$, for all
$j\in\fso{n-1}$, and $z\nleq W_{j,n}(\vec b,\vec b')$, for all
$j\in\fso{n-2}$, and
$z\nleq(a_1\vee b'_1)\wedge(a'_1\wedge b_1)$, say,
$z\nleq a_1\vee b'_1$. By Lemma~\ref{L:JIT2Id1}, there are
$x_1\leq a^*_1$, \dots, $x_m\leq a^*_m$, $x'_1\leq a'_1$, \dots,
$x'_m\leq a'_m$, $y_1\leq b^*_1$, \dots, $y_n\leq b^*_n$,
$y'_1\leq b'_1$, \dots, $y'_n\leq b'_n$ in $\Sigma$ such that,
putting $x_0=y_0=z$, both pairs
 \begin{align*}
 \sigma&=(\famm{x_i}{0\leq i\leq m},\famm{x'_i}{1\leq i\leq m}),\\
 \tau&=(\famm{y_j}{0\leq j\leq n},\famm{y'_j}{1\leq j\leq n})
 \end{align*}
are Stirlitz tracks. By assumption, the pair $(\sigma,\tau)$ is not
a bi-Stirlitz track, whence $z\nleq x_1\vee y_1$. Furthermore, from
$z\nleq a_1\vee b'_1$ it follows that $z\nleq x_1\vee y'_1$ (observe
that $x_1\leq a^*_1\leq a_1$). However, from the fact that 
$z\leq x_1\vee x'_1,y_1\vee y'_1$ are nontrivial join-covers and
\Boj\ it follows that either $z\leq x_1\vee y_1$ or
$z\leq x_1\vee y'_1$, a contradiction.
\end{proof}

\begin{corollary}\label{C:MainHmnP}
Let $m$ and $n$ be positive integers, let $P$ be a poset. Then
$\Co(P)$ satisfies \Ht{m,n} if{f} $\lh P\leq m+n-1$.
\end{corollary}

\begin{proof}
Suppose first that $P$ contains a $m+n+1$-element chain, say,
 \[
 x_m\tr\cdots\tr x_1\tr x_0=y_0\tr y_1\tr\cdots\tr y_n.
 \]
Then both pairs $\sigma$ and $\tau$ defined as
 \begin{align*}
 \sigma&=(\famm{\set{x_i}}{0\leq i\leq m},
 \famm{\set{y_1}}{1\leq i\leq m})\\
 \tau&=(\famm{\set{y_j}}{0\leq j\leq n},
 \famm{\set{x_1}}{1\leq j\leq n})
 \end{align*}
are Stirlitz tracks with the same base
$\set{x_0}=\set{y_0}\leq\set{x_1}\vee\set{y_1}$, hence
$(\sigma,\tau)$ is a bi-Stirlitz track of index $(m,n)$. By
Proposition~\ref{P:MainHmn}, $\Co(P)$ does not satisfy~\Ht{m,n}.

Conversely, suppose that $P$ does not contain any $m+n+1$-element
chain. By Proposition~\ref{P:MainHmn}, in order to prove that
$\Co(P)$ satisfies \Ht{m,n}, it suffices to prove that it has no
bi-Stirlitz track of index $(m,n)$ with entries in
$\Sigma=\setm{\set{p}}{p\in P}$. Let
 \begin{align*}
 \sigma&=(\famm{\set{x_i}}{0\leq i\leq m},
 \famm{\set{x'_i}}{1\leq i\leq m})\\
 \tau&=(\famm{\set{y_j}}{0\leq j\leq n},
 \famm{\set{y'_j}}{1\leq j\leq n})
 \end{align*}
be pairs such that $(\sigma,\tau)$ is such a bi-Stirlitz track. By
an argument similar as the one used in the proof of
Corollary~\ref{C:MainHnP}, since $\sigma$ is a Stirlitz track,
either $x'_1\tr x_0\tr\cdots\tr x_m$ or
$x_m\tr\cdots\tr x_0\tr x'_1$; without loss of generality,
the second possibility occurs. Similarly, since $\tau$ is a
Stirlitz track, either $y'_1\tr y_0\tr\cdots\tr y_n$ or
$y_n\tr\cdots\tr y_0\tr y'_1$. If the second possibility
occurs, then $y_1\tr y_0=x_0$ and $x_1\tr x_0$ while
$\set{x_0}\leq\set{x_1}\vee\set{y_1}$, a contradiction. Therefore,
the first possibility occurs, hence
 \[
 x_m\tr\cdots\tr x_1\tr x_0=y_0\tr y_1\tr\cdots\tr y_n
 \]
is a $m+n+1$-element chain in $P$, a contradiction.
\end{proof}

Now let us recall some results of \cite{SeWe1}. In case $L$ belongs
to the variety $\SUB$, so does the lattice $\hL=\Fil L$ of all
filters of $L$ partially ordered by reverse inclusion
(see Section~\ref{S:Letna}), and $\J(\hL)$ is a join-seed of~$\hL$.
Furthermore, one can construct two posets $R$ and $\Gamma$ with the
following properties:

\begin{enumerate}
\item There are natural embeddings $\varphi\colon L\into\Co(R)$ and
$\psi\colon L\into\Co(\Gamma)$, and they preserve the existing
bounds.

\item $R$ is finite in case $L$ is finite.

\item $\Gamma$ is tree-like (as defined in Section~\ref{S:Basic},
see also \cite{SeWe1}).

\item There exists a natural map $\pi\colon\Gamma\to R$ such that
$\alpha\prec\beta$ in $\Gamma$ implies that
$\pi(\alpha)\prec\pi(\beta)$ in $R$. In particular, $\pi$ is
order-preserving.

\item $\psi(x)=\pi^{-1}[\varphi(x)]$, for all $x\in L$.
\end{enumerate}

The main theorem of this section is the following:

\begin{theorem}\label{T:SUBnVar}
Let $n$ be a positive integer,
let $L$ be a lattice that belongs to the variety $\SUB$. Consider
the posets $R$ and $\Gamma$ constructed in \textup{\cite{SeWe1}} from
$\hL$. Then the following are equivalent:
\begin{enumerate}
\item $\lh R\leq n$;

\item $\lh\Gamma\leq n$;

\item there exists a poset $P$ such that $\lh P\leq n$ and $L$
embeds into $\Co(P)$;

\item $L$ satisfies the identities \Ht{n} and \Ht{k,n+1-k} for
$1<k<n$;

\item $L$ satisfies the identities \Ht{n} and \Ht{k,n+1-k} for
$1\leq k\leq n$.

\end{enumerate}
\end{theorem}

\begin{proof}
(i)$\Rightarrow$(ii) Suppose that $\lh R\leq n$, we prove that
$\lh\Gamma\leq n$. Otherwise, there exists a $n+2$-element chain
$\alpha_0\prec\cdots\prec\alpha_{n+1}$ in $\Gamma$, thus, applying
the map $\pi$, we obtain a $n+2$-element chain
$\pi(\alpha_0)\prec\cdots\prec\pi(\alpha_{n+1})$ in $R$, a
contradiction.

(ii)$\Rightarrow$(iii) Since $L$ embeds into $\Co(\Gamma)$, it
suffices to take $P=\Gamma$.

(iii)$\Rightarrow$(iv) follows immediately from
Corollaries~\ref{C:MainHnP} and \ref{C:MainHmnP}.

(iv)$\Rightarrow$(v) Suppose that $L$ satisfies the identities
\Ht{n} and \Ht{k,n+1-k} for $1<\nobreak k<\nobreak n$; then so does
the filter lattice $\hL$ of $L$. Since $\hL$ satisfies \Ht{n}, it
has no Stirlitz track of length $n$ (see Proposition~\ref{P:MainHn}),
thus, \emph{a fortiori}, it has no bi-Stirlitz track of index either
$(n,1)$ or $(1,n)$. Since $\J(\hL)$ is a join-seed of $\hL$, it
follows from Proposition~\ref{P:MainHmn} that $\hL$ satisfies
both \Ht{n,1} and \Ht{1,n}.

(v)$\Rightarrow$(i) Suppose that $L$ satisfies the identities
\Ht{n} and \Ht{k,n+1-k} for $1\leq\nobreak k\leq\nobreak n$; then so
does the filter lattice $\hL$ of $L$. We prove that $\lh R\leq n$.
Otherwise, $R$ has an oriented path $\fr=\seq{r_0,\dots,r_{n+1}}$ of
length $n+2$, that is, $r_i\prec r_{i+1}$, for all $i\in\fso{n}$.
By \cite[Lemma~6.4]{SeWe1}, we can assume that
$\fr$ is `reduced'. If there are $n$ successive values of the $r_i$
that are of the form $\seq{a_i,b_i,\eps}$ for a constant
$\eps\in\set{+,-}$, then, by
\cite[Lemma~6.1]{SeWe1}, there exists a Stirlitz track of length $n$
in $\hL$ (with entries in $\J(\hL)$), which contradicts the
assumption that $\hL$ satisfies \Ht{n} and
Proposition~\ref{P:MainHn}. Therefore, $\fr$ has the form
 \[
 \seq{\seq{a_{k-1},a_k,-},\dots,\seq{a_0,a_1,-},\seq{p},
 \seq{b_0,b_1,+},\dots,\seq{b_{l-1},b_l,+}}
 \]
for some positive integers $k$ and $l$ and elements $a_0$,
\dots, $a_k$, $b_0$, \dots, $b_l$ of $\J(\hL)$. By
\cite[Lemma~6.1]{SeWe1}, there are Stirlitz tracks of the form
 \begin{align*}
 \sigma&=(\famm{a_i}{0\leq i\leq k},\famm{a'_i}{1\leq i\leq k}),\\
 \tau&=(\famm{b_j}{0\leq j\leq l},\famm{b'_j}{1\leq j\leq l})
 \end{align*}
for elements $a'_1$, \dots, $a'_k$, $b'_1$, \dots, $b'_l$ of
$\J(\hL)$. Observe that $p=a_0=b_0$. Furthermore, from
$\seq{a_0,a_1,-}\prec\seq{p}\prec\seq{b_0,b_1,+}$ and the
definition of $\prec$ on $R$ it follows that $p\leq a_1\vee b_1$.
Therefore, $(\sigma,\tau)$ is a bi-Stirlitz track of index $(k,l)$
with $k+l=n+1$ in $\hL$, which contradicts the assumption that $\hL$
satisfies \Ht{k,l} and Proposition~\ref{P:MainHmn}.
\end{proof}

The main result of \cite{SeWe1} is that $\SUB$ is a finitely based
variety of lattices. We thus obtain the following:

\begin{corollary}\label{C:SUBnVar1}
Let $n$ be a positive integer.
The class $\SUB_n$ of all lattices $L$ that can be embedded into
$\Co(P)$ for a poset $P$ of length at most $n$ is a finitely based
variety, defined by the identities \St, \Ud, \Bo, \Ht{n}, and
\Ht{k,n+1-k} for $1<k<n$.
\end{corollary}

Since finiteness of $L$ implies finiteness of $R$, we also obtain
the following:

\begin{corollary}\label{C:SUBnVar2}
Let $n$ be a positive integer.
A finite lattice $L$ belongs to $\SUB_n$ if{f} it can be embedded
into $\Co(P)$ for some finite poset $P$ of length at most $n$.
\end{corollary}

For a positive integer $m$, denote by $\boldsymbol{m}$ the
$m$-element chain. As a consequence of Corollaries~\ref{C:MainHnP}
and \ref{C:MainHmnP} and of Theorem~\ref{T:SUBnVar}, we obtain
immediately the following:

\begin{corollary}\label{C:chainsHn}
For positive integers $m$ and $n$, $\Co(\boldsymbol{m})$ belongs to
$\SUB_n$ if{f}\linebreak
$m\leq n+1$. In particular, $\SUB_n$ is a
proper subvariety of $\SUB_{n+1}$, for every positive integer $n$.
\end{corollary}

\section{Non-local finiteness of $\SUB_3$}\label{S:NonLF3}

We have seen in Section~\ref{S:LocFin2} that the variety $\SUB_2$ is
locally finite. In contrast with this, we shall now prove the
following:

\begin{theorem}\label{T:NonLF3}
There exists an infinite, three-generated lattice in $\SUB_3$. Hence
$\SUB_n$ is not locally finite for $n\geq3$.
\end{theorem}

\begin{proof}
Let $P$ be the poset diagrammed on Figure~2.

\begin{figure}[htb]
\includegraphics{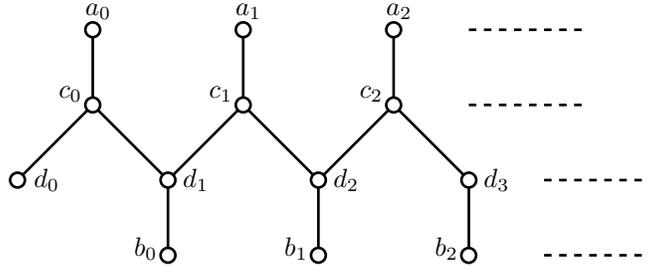}
\caption{An infinite poset of length $3$}
\end{figure}

We observe that the length of $P$ is $3$. We define
order-convex subsets $A$, $B$, $C$ of $P$ as follows:
 \[
 A=\setm{a_n}{n<\omega},\quad B=\set{d_0}\cup\setm{b_n}{n<\omega},
 \quad C=\setm{c_n}{n<\omega}\cup\setm{d_n}{n<\omega}.
 \]
We put $A_0=A$, $B_0=B$, $A_{n+1}=A\vee(B_n\cap C)$, and
$B_{n+1}=B\vee(A_n\cap C)$, for all $n<\omega$. A straightforward
computation yields that both $c_n$ and $d_n$ belong to
$A_{2n+1}\setminus A_{2n}$, for all $n<\omega$. Hence the sublattice
of $\Co(P)$ generated by $\set{A,B,C}$ is infinite.
\end{proof}

\section{Open problems}\label{S:Pbs}

So far we have studied the following $(\omega+1)$-chain of
varieties:
 \begin{equation}\label{Eq:VarChain}
 \mathbf{D}=\SUB_1\subset\SUB_2\subset\SUB_3\subset\cdots
 \subset\SUB_n\subset\cdots\subset\SUB.
 \end{equation}
We do not know the answer to the following simple question, see also
Problem~1 in~\cite{SeWe1}:

\begin{problem}\label{Pb:joinSUBn}
Is $\SUB$ the quasivariety join of all the $\SUB_n$,
for $n>0$?
\end{problem}

Every variety from the chain \eqref{Eq:VarChain} is the
variety $\SUB(\mathcal{K})$ generated by all $\Co(P)$, where
$P\in\mathcal{K}$, for some class $\mathcal{K}$ of posets.

\begin{problem}\label{Pb:SpKVar}
Can one classify all the varieties of the form $\SUB(\mathcal{K})$?
In particular, are there only countably many such varieties?
\end{problem}

\begin{problem}\label{Pb:CpleteSubCoP}
What are the \emph{complete} sublattices of the lattices of the form
$\Co(P)$ for some poset $P$?
\end{problem}

\begin{problem}\label{Pb:FreeSUB2}
Give an estimate for the cardinality of the free lattice in
$\SUB_2$ on $m$ generators, for a positive integer $m$.
\end{problem}

\begin{problem}\label{Pb:SubvSUB2}
Classify all the subvarieties of $\SUB_2$.
\end{problem}

\section{Acknowledgments}
This work was partially completed while both authors were visiting
the Charles University during the fall of 2001. Excellent conditions
provided by the Department of Algebra are highly appreciated.
Special thanks are due to Ji\v{r}\'\i\ T\r{u}ma and
Vaclav Slav\'\i k.

\end{document}